\documentclass[12pt,reqno]{amsart}

\usepackage{mathdots}

\usepackage{color}

\usepackage{pdfsync}

\usepackage{enumitem}

\usepackage{wasysym}

\usepackage{amssymb}

\usepackage{mathrsfs}

\DeclareMathAlphabet{\mathpzc}{OT1}{pzc}{m}{it}

\usepackage[all]{xy}

\usepackage{tikz}
\usetikzlibrary{arrows,decorations.pathmorphing,decorations.pathreplacing,positioning,shapes.geometric,shapes.misc,decorations.markings,decorations.fractals,calc,patterns}

\usepackage{float}

\usepackage[bottom]{footmisc}

\usepackage{mathtools}

\entrymodifiers={+!!<0pt,\fontdimen22\textfont2>}

\def\cA{\mathscr{A}}
\def\cB{\mathscr{B}}
\def\cC{\mathscr{C}}
\def\cD{\mathscr{D}}

\def\cH{\mathscr{H}}

\def\cK{\mathscr{K}}

\def\cM{\mathscr{M}}

\def\cP{\mathscr{P}}

\def\cS{\mathscr{S}}
\def\cT{\mathscr{T}}

\def\cX{\mathscr{X}}
\def\cY{\mathscr{Y}}

\def\BC{\mathbb{C}}

\def\BZ{\mathbb{Z}}

\def\add{\operatorname{add}}

\def\adots{\mathinner{\mkern1mu\raise1.0pt\vbox{\kern7.0pt\hbox{.}}\mkern2mu\raise4.0pt\hbox{.}\mkern2mu\raise7.0pt\hbox{.}\mkern1mu}}

\def\ast{{\textstyle *}}

\def\Coker{\operatorname{Coker}}

\def\dddots{\mathinner{\mkern1mu\raise10.0pt\vbox{\kern7.0pt\hbox{.}}\mkern2mu\raise5.3pt\hbox{.}\mkern2mu\raise1.0pt\hbox{.}\mkern1mu}}
\def\dddotssmall{\mathinner{\mkern1mu\raise7.0pt\vbox{\kern7.0pt\hbox{.}}\mkern-1mu\raise4pt\hbox{.}\mkern-1mu\raise1.0pt\hbox{.}\mkern1mu}}

\def\deg{\operatorname{deg}}

\def\Db{\cD^{\operatorname{b}}}

\def\H{\operatorname{H}}

\def\Hom{\operatorname{Hom}}

\def\id{\operatorname{id}}

\def\Kb{\cK^{\operatorname{b}}}

\def\Ker{\operatorname{Ker}}

\def\mod{\operatorname{mod}}
\def\Mod{\operatorname{Mod}}

\def\prj{\operatorname{prj}}
\def\Prj{\operatorname{Prj}}

\def\SL2{\operatorname{SL}_2}

\numberwithin{equation}{section}

\newtheorem{Lemma}{Lemma}[section]
\newtheorem{Theorem}[Lemma]{Theorem}
\newtheorem{Proposition}[Lemma]{Proposition}

\theoremstyle{definition}
\newtheorem{Definition}[Lemma]{Definition}

\newtheorem{Remark}[Lemma]{Remark}

\newtheorem{Example}[Lemma]{Example}

\newtheorem*{bfhpg*}{}

\begin{document}

\setlength{\parindent}{0pt}
\setlength{\parskip}{7pt}

\title[Co-t-structures decade]{Co-t-structures: The First Decade}

\author{Peter J\o rgensen}
\address{School of Mathematics and Statistics,
Newcastle University, Newcastle upon Tyne NE1 7RU, United Kingdom}
\email{peter.jorgensen@ncl.ac.uk}
\urladdr{http://www.staff.ncl.ac.uk/peter.jorgensen}




\begin{abstract} 

  Co-t-structures were introduced about ten years ago as a type of
  mirror image of t-structures.  Like t-structures, they permit to
  divide an object in a triangulated category $\cT$ into a ``left
  part'' and a ``right part'', but there are crucial differences.  For
  instance, a bounded t-structure gives rise to an abelian subcategory
  of $\cT$, while a bounded co-t-structure gives rise to a so-called
  silting subcategory.

  This brief survey will emphasise three philosophical points.  First,
  bounded t-structures are akin to the canonical example of ``soft''
  truncation of complexes in the derived category.  Secondly, bounded
  co-t-structures are akin to the canonical example of ``hard''
  truncation of complexes in the homotopy category.

  Thirdly, a triangulated category $\cT$ may be skewed towards
  t-structures or co-t-structures, in the sense that one type of
  structure is more useful than the other for studying $\cT$.  In
  particular, we think of derived categories as skewed towards
  t-structures, and of homotopy categories as skewed towards
  co-t-structures.

\end{abstract}

\maketitle

\setcounter{section}{-1}
\section{Introduction}
\label{sec:introduction}

The notion of co-t-structure in a triangulated category $\cT$ was
introduced independently by Bondarko and Pauksztello, see Definition
\ref{def:co-t-structure}.  It is a mirror image of the classic notion of
t-structure due to Beilinson, Bernstein, and Deligne, see Definition
\ref{def:t-structure}.

Given an object $t \in \cT$, both types of structure give a way to
divide $t$ into a ``left part'' and a ``right part''.  This is
exemplified by dividing a complex of modules into a left and a right
part by ``soft'' or ``hard'' truncation, see Figures \ref{fig:1} and
\ref{fig:2}.

Each case gives a triangle $u \rightarrow t \rightarrow v$.
Crucially, for a t-structure, $u$ is the left part of $t$ and $v$ the
right part; for a co-t-structure, vice versa.  This reversal leads to
a number of differences, and the theories of t-structures and
co-t-structures are far from being simple mirrors of each other.  For
instance, while a bounded t-structure induces an abelian subcategory
of $\cT$, a bounded co-t-structure induces a so-called silting
subcategory; see Definition \ref{def:silting}.

In the decade since their inception, the theory of co-t-structures has
grown considerably.  This brief survey is far from encyclopedic, but
has the important goal of communicating three philosophical points:
\begin{enumerate}
\setlength\itemsep{4pt}
\setcounter{enumi}{0}

  \item  Bounded t-structures are akin to soft truncation in the bounded
  derived category of an abelian category.

  \item  Bounded co-t-structures are akin to hard truncation in the
  bounded homotopy ca\-te\-go\-ry of an additive category.

  \item  A triangulated category $\cT$ may be skewed in the direction of
  t-structures or co-t-struc\-tu\-res, in the sense that one type of
  structure is more useful than the other for studying $\cT$.

\end{enumerate}
We explain these points in the next three subsections.  Throughout,
$\cT$ is a triangulated category with $\Hom$-spaces $\cT( -,- )$ and
suspension functor $\Sigma$.  If $\cM$ is an abelian category then
$\Db( \cM )$ is the derived category of bounded complexes over $\cM$,
and if $\cP$ is an additive category then $\Kb( \cP )$ is the homotopy
category of bounded complexes over $\cP$.

{\bf (i) Bounded t-structures are akin to soft truncation in the
  bounded derived category.}  Let $R$ be a ring, $\cM = \Mod R$ the
category of left modules over $R$.  Each complex $t \in \Db( \cM )$
has soft truncations $a$ and $b$ as shown in Figure \ref{fig:1}, and
there is a triangle $a \rightarrow t \rightarrow b$ in $\Db( \cM )$.
\begin{figure}
\[
  \xymatrix 
  {
    a \;=\; \cdots \ar@<-3.5ex>[d] \ar[r] & M^{-2} \ar[r] \ar[d] & M^{-1} \ar[r] \ar[d] & \Ker d^0 \ar[r] \ar[d] & 0 \ar[r] \ar[d] & 0 \ar[r] \ar[d] & 0 \ar[r] \ar[d] & \cdots \\
    t \;=\; \cdots \ar@<-3.5ex>[d] \ar[r] & M^{-2} \ar[r] \ar[d] & M^{-1} \ar[r] \ar[d] & M^0 \ar^{d^0}[r] \ar[d] & M^1 \ar[r] \ar[d] & M^2 \ar[r] \ar[d] & M^3 \ar[r] \ar[d] & \cdots \\
    b \;=\; \cdots \ar[r] & 0 \ar[r] & 0 \ar[r] & 0 \ar[r] & \Coker d^0 \ar[r] & M^2 \ar[r] & M^3 \ar[r] & \cdots
  }
\]
  \caption{A complex $t = \cdots \rightarrow M^{-1} \rightarrow M^0
  \stackrel{d^0}{\rightarrow} M^1 \rightarrow M^2 \rightarrow \cdots$
  has soft truncations $a$ and $b$.  There is a triangle $a
  \rightarrow t \rightarrow b$ in the derived category.  Each vertical
  module homomorphism is either the identity, a canonical inclusion or
  surjection, or zero.}
\label{fig:1}
\end{figure}
The full subcategories $\cA$ and $\cB$ consisting of complexes
isomorphic to such truncations satisfy Definition
\ref{def:t-structure} and hence form a bounded t-structure, sometimes
known as the standard t-structure, see Example \ref{exa:t-structure}.

Up to isomorphism, the heart $\cH = \cA \cap \Sigma \cB$ consists of
complexes concentrated in degree $0$.  The heart is an abelian
subcategory of $\Db( \cM )$ which is equivalent to $\cM$.  Each $t \in
\Db( \cM )$ permits a ``tower'' as shown in Proposition
\ref{pro:t-tower} due to Beilinson, Bernstein, and Deligne.  This
expresses how to build $t$ from objects of the form $\Sigma^i h$ with
$h \in \cH$.  The objects can be taken to be $\Sigma^i\H^{ -i }( t )$
and then the tower shows how $t$ is built from its cohomology modules,
see Example \ref{exa:t-tower}.

In general, if $( \cA,\cB )$ is a bounded t-structure in a
triangulated category $\cT$, then each $t \in \cT$ still permits a
triangle $a \rightarrow t \rightarrow b$ with $a \in \cA$, $b \in
\cB$.  The heart $\cH$ is still an abelian subcategory of $\cT$ by
Theorem \ref{thm:abelian} due to Beilinson, Bernstein, and Deligne,
and each $t \in \cT$ still has the tower in Proposition
\ref{pro:t-tower}.

Working in such a setup is akin to working with soft truncation in
$\Db( \cM )$.

{\bf (ii) Bounded co-t-structures are akin to hard truncation in the
  bounded homotopy category.}  Let $R$ be a ring, $\cP = \Prj R$ the
category of projective left modules over $R$.  Each complex $t \in
\Kb( \cP )$ has hard truncations $x$ and $y$ as shown in Figure
\ref{fig:2}, and there is a triangle $x \rightarrow t \rightarrow y$
in $\Kb( \cP )$.
\begin{figure}
\[
  \xymatrix 
  {
    x \;=\; \cdots \ar@<-3.5ex>[d] \ar[r] & 0 \ar[r] \ar[d] & 0 \ar[r] \ar[d] & 0 \ar[r] \ar[d] & P^0 \ar[r] \ar[d] & P^1 \ar[r] \ar[d] & P^2 \ar[r] \ar[d] & \cdots \\
    t \;=\; \cdots \ar@<-3.5ex>[d] \ar[r] & P^{-3} \ar[r] \ar[d] & P^{-2} \ar[r] \ar[d] & P^{-1} \ar[r] \ar[d] & P^0 \ar[r] \ar[d] & P^1 \ar[r] \ar[d] & P^2 \ar[r] \ar[d] & \cdots \\
    y \;=\; \cdots \ar[r] & P^{-3} \ar[r] & P^{-2} \ar[r] & P^{-1} \ar[r] & 0 \ar[r] & 0 \ar[r] & 0 \ar[r] & \cdots
  }
\]
  \caption{A complex $t = \cdots \rightarrow P^{-2} \rightarrow P^{-1}
  \rightarrow P^0 \rightarrow P^1 \rightarrow \cdots$ has hard
  truncations $x$ and $y$.  There is a triangle $x \rightarrow t
  \rightarrow y$ in the homotopy category.  Each vertical module
  homomorphism is the identity or zero.}
\label{fig:2}
\end{figure}
The full subcategories $\cX$ and $\cY$ consisting of complexes
isomorphic to such truncations satisfy Definition
\ref{def:co-t-structure} and hence form a bounded co-t-structure,
sometimes known as the standard co-t-structure, see Example
\ref{exa:co-t-structure}.

Up to isomorphism, the co-heart $\cC = \cX \cap \Sigma^{-1} \cY$
consists of complexes concentrated in degree $0$.  The co-heart is an
additive subcategory of $\Kb( \cP )$ which is equivalent to $\cP$.

The co-heart is not in general abelian, but it does have the strong
property of being a so-called silting subcategory of $\Kb( \cP )$, see
Definition \ref{def:silting} and Theorem \ref{thm:MSSS} due to Mendoza
Hern\'{a}ndez et.al.  Such subcategories are structurally important.

Each $t = \cdots \rightarrow P^{-3} \rightarrow P^{-2} \rightarrow
P^{-1} \rightarrow P^0 \rightarrow P^1 \rightarrow P^2 \rightarrow
\cdots$ in $\Kb( \cP )$ permits a ``tower'' as shown in Proposition
\ref{pro:co-t-tower} due to Bondarko.  This expresses how to build $t$
from objects of the form $\Sigma^i c$ with $c \in \cC$.  The objects
can be taken to be $\Sigma^iP^{ -i }$ and then the tower shows how $t$
is built from its constituent modules, see Example
\ref{exa:co-t-tower}.

In general, if $( \cX,\cY )$ is a bounded co-t-structure in a
triangulated category $\cT$, then each $t \in \cT$ still permits a
triangle $x \rightarrow t \rightarrow y$ with $x \in \cX$, $y \in
\cY$.  The co-heart $\cC$ is still a silting subcategory of $\cT$ by
Theorem \ref{thm:MSSS}, and each $t \in \cT$ still has the tower in
Proposition \ref{pro:co-t-tower}.

Working in such a setup is akin to working with hard truncation in
$\Kb( \cP )$.

{\bf (iii) Triangulated categories skewed in the direction of
  t-structures or co-t-structures.}
When studying a given triangulated category, bounded co-t-structures
may be more useful than bounded t-structures, simply because there are
none of the latter.  See Section \ref{sec:skewed} for a toy example. 

Section \ref{sec:KY} shows a subtler skewing phenomenon.  Let
$\Lambda$ be a finite dimensional $\BC$-algebra, $\cM = \mod\,
\Lambda$ the category of finite dimensional left modules over
$\Lambda$, and $\cP = \prj\, \Lambda$ the category of finite
dimensional projective left modules over $\Lambda$.  Theorem
\ref{thm:KY}, due to K\"{o}nig and Yang, shows a bijection between all
bounded co-t-structures in $\Kb( \cP )$ and the bounded t-structures
in $\Db( \cM )$ whose hearts are length categories (these are, in a
sense, the ``algebraic'' t-structures; for instance, they permit a
nice mutation theory).

In itself, this does not imply that either category has more
t-structures than co-t-structures or vice versa, but we think of it as
indicating that when going from derived to homotopy categories, the
role played by t-structures is taken over by co-t-structures.

{\bf Summing up, } Sections \ref{sec:t-structures} and
\ref{sec:co-t-structures} show the definitions of t-structures and
co-t-structures, emphasising the similarities with soft and hard
truncation of complexes of modules.  Section \ref{sec:towers} shows
the ``towers'' whereby general objects can be built from objects in
the (co-)heart.  Sections \ref{sec:skewed} and \ref{sec:KY} illustrate
how a triangulated category may be skewed towards t-structures or
co-t-structures.

Section \ref{sec:mutation} shows the silting mutation of Aihara and
Iyama which mutates between bounded co-t-structures.  It permits the
definition of the so-called silting quiver which is a combinatorial
picture of how the bounded co-t-structures fit inside a fixed
triangulated category.

\section{t-structures}
\label{sec:t-structures}

The following definition was made in \cite[def.\ 1.3.1]{BBD}.

\begin{Definition}
[Beilinson, Bernstein, and Deligne]
\label{def:t-structure}
A {\em t-structure} in the triangulated category $\cT$ is a pair $(
\cA,\cB )$ of full subcategories, closed under isomorphisms, direct
sums, and direct summands, which satisfy the following conditions.
\begin{enumerate}
\setlength\itemsep{4pt}
\setcounter{enumi}{0}

  \item  $\Sigma\cA \subseteq \cA$ and $\Sigma^{ -1 }\cB \subseteq \cB$.

  \item  $\cT( \cA,\cB ) = 0$.

  \item  For each object $t \in \cT$ there is a triangle
  $a \rightarrow t \rightarrow b$ with $a \in \cA$, $b \in \cB$. 

\end{enumerate}
The {\em heart} is $\cH = \cA \cap \Sigma\cB$.

The t-structure is called {\em bounded} if
\[
\tag*{$\Box$}
  \bigcup_{i \in \BZ} \Sigma^i \cA = \bigcup_{i \in \BZ} \Sigma^i \cB = \cT.
\]
\end{Definition}

The objects $a$ and $b$ in Definition \ref{def:t-structure}(iii)
depend functorially on $t$.  The resulting functor $t \mapsto a$ is a
right-adjoint to the inclusion $\cA \hookrightarrow \cT$.  Similarly,
$t \mapsto b$ is a left adjoint to the inclusion $\cB \hookrightarrow
\cT$.  See \cite[prop.\ 1.3.3]{BBD}.

The following is the canonical example of a t-structure.

\begin{Example}
[The standard t-structure]
\label{exa:t-structure}
Let $R$ be a ring, $\cM = \Mod\, R$ the category of left modules over
$R$.  Let $\cA$ and $\cB$ be the isomorphism closures in the bounded
derived category $\Db( \cM )$ of the subsets
\[
  \xymatrix @-2pc @C=0.7pc @R=0.3pc
  {
    \{\, \cdots \ar[r] & M^{-2} \ar[r] & M^{-1} \ar[r] & M^0 \ar[r] & 0 \ar[r] & 0 \ar[r] & 0 \ar[r] & \cdots \,|\, M^i \in \cM, \, M^i = 0 \mbox{ for } i \ll 0 \,\} \lefteqn{,} \\
    \{\, \cdots \ar[r] & 0 \ar[r] & 0 \ar[r] & 0 \ar[r] & M^1 \ar[r] & M^2 \ar[r] & M^3 \ar[r] & \cdots \,|\, M^i \in \cM, \, M^i = 0 \mbox{ for } i \gg 0 \,\} \lefteqn{.} 
  }
\]
We will show that if $\cA$ and $\cB$ are viewed as full subcategories,
then $( \cA,\cB )$ is a bounded t-structure with heart $\cH$
equivalent to $\cM$.

It is easy to show
\begin{equation}
\label{equ:t-structure_by_H}
\vcenter{
  \xymatrix @R=0.23pc
  {
    \cA = \{\, M \in \Db( \cM ) \,|\, \H^i( M ) = 0 \mbox{ for } i \geqslant 1 \,\}, \\
    \cB = \{\, M \in \Db( \cM ) \,|\, \H^i( M ) = 0 \mbox{ for } i \leqslant 0 \,\}.
  }
        }
\end{equation}
This description clearly implies that $\cA$ and $\cB$ are closed under
isomorphisms, direct sums, and direct summands.

We next check the conditions in Definition \ref{def:t-structure}.
Condition (i) is immediate.  Condition (ii) requires that if objects
\[
  \xymatrix @-2pc @C=0.7pc @R=0.3pc
  {
    a = \cdots \ar[r] & M^{-2} \ar[r] & M^{-1} \ar[r] & M^0 \ar[r] & 0 \ar[r] & 0 \ar[r] & 0 \ar[r] & \cdots \lefteqn{,} \\
    b = \cdots \ar[r] & 0 \ar[r] & 0 \ar[r] & 0 \ar[r] & M^1 \ar[r] & M^2 \ar[r] & M^3 \ar[r] & \cdots
  }
\]
in $\cA$ and $\cB$ are given, then $\Hom_{ \Db( \cM ) }( a,b ) = 0$.
This can be shown by noting that $a$ has a projective resolution
\[
  p = \cdots \rightarrow P^{-2} \rightarrow P^{-1} \rightarrow P^0 \rightarrow 0 \rightarrow 0 \rightarrow 0 \rightarrow \cdots,
\]
and that
\[
  \Hom_{ \Db( \cM ) }( a,b )
  = \Hom_{ \cK( \cM ) }( p,b )
  = 0.
\]
Here $\cK( \cM )$ is the homotopy category of complexes over $\cM$,
and the second $=$ holds because in each degree, either the complex
$p$ or the complex $b$ is zero.  The triangle in Definition
\ref{def:t-structure}(iii) can be obtained by soft truncation of
the object $t = \cdots \rightarrow M^{-1} \rightarrow M^0
\stackrel{d^0}{\rightarrow} M^1 \rightarrow M^2 \rightarrow \cdots$ in
$\Db( \cM )$, see Figure \ref{fig:1} in the introduction.

It is immediate that the t-structure $( \cA,\cB )$ is bounded.

Finally, it follows from Equation \eqref{equ:t-structure_by_H} that the heart
$\cH = \cA \cap \Sigma\cB$ is
\[
  \cH = \{\, M \in \Db( \cM ) \,|\, \H^i( M ) = 0 \mbox{ for } i \neq 0 \,\},
\]
and this subcategory of $\Db( \cM )$ is equivalent to $\cM$. 
\hfill $\Box$
\end{Example}

The following pivotal result is one of the motivations for the
definition of t-structures.  It was proved in \cite[thm.\ 1.3.6]{BBD}.

\begin{Theorem}
[Beilinson, Bernstein, and Deligne]
\label{thm:abelian}
Let $( \cA,\cB )$ be a t-structure in $\cT$.  Then the heart $\cH =
\cA \cap \Sigma \cB$ is an abelian subcategory of $\cT$.
\end{Theorem}

\begin{Example}
\label{exa:little_t-structure}
Let $\Lambda = \BC A_2$ be the path algebra of the quiver
\begin{equation}
\label{equ:A2}
  A_2 \;\;=\;\; 1 \rightarrow 2
\end{equation}
and let $\cM = \mod\,\Lambda$ be the category of finite dimensional
left modules over $\Lambda$.

There is a bounded t-structure $( \cA,\cB )$ in $\Db( \cM )$ defined
by Equation \eqref{equ:t-structure_by_H}; this is shown by the same
method as in Example \ref{exa:t-structure}.

There are three isomorphism classes of indecomposable objects in $\cM$
given by the following representations of the quiver $A_2$.
\[
  0 \rightarrow \BC
  \;\;,\;\;
  \BC \stackrel{\id}{\rightarrow} \BC
  \;\;,\;\;
  \BC \rightarrow 0
\]
They induce isomorphism classes $x_0$, $x_1$, $x_2$ of indecomposable
objects in $\Db( \cM )$, and we define further isomorphism classes
recursively by $\Sigma x_i = x_{ i+3 }$ for $i \in \BZ$.

The Auslander--Reiten quiver of $\Db( \cM )$ looks as follows, where
red and green vertices show $\cA$ and $\cB$.
\begin{equation}
\label{equ:AR_quiver}
\vcenter{
  \xymatrix @-1.0pc @! {
    & {\color{green} x_{ -3 }} \ar[dr] & & {\color{green} x_{ -1 }} \ar[dr] & & {\color{red} x_1} \ar[dr] & & {\color{red} x_3} \ar[dr] & & {\color{red} x_5} \ar[dr] & \\
    \cdots \ar[ur] & & {\color{green} x_{ -2 }} \ar[ur] & & {\color{red} x_0} \ar[ur] & & {\color{red} x_2} \ar[ur] & & {\color{red} x_4} \ar[ur] & & \cdots
               }
        }
\end{equation}
Note that if $t \in \Db( \cM )$ is indecomposable, then $t$ is in
$\cA$ or in $\cB$, so the triangle in Definition
\ref{def:t-structure}(iii) is trivial in the sense that it reads $t
\rightarrow t \rightarrow 0$ or $0 \rightarrow t \rightarrow t$.  

The heart $\cH = \cA \cap \Sigma \cB$ is determined by
\[
\tag*{$\Box$}
  \cH = \add( x_0 \oplus x_1 \oplus x_2 ).
\]
\end{Example}

\section{Co-t-structures}
\label{sec:co-t-structures}

The following definition was made in \cite[def.\ 1.1.1]{B} and
\cite[def.\ 2.4]{P}.

\begin{Definition}
[Bondarko and Pauksztello]
\label{def:co-t-structure}
A {\em co-t-structure} in $\cT$ is a pair $( \cX,\cY )$ of full
subcategories, closed under isomorphisms, direct sums, and
direct summands, which satisfy the following conditions.
\begin{enumerate}
\setlength\itemsep{4pt}
\setcounter{enumi}{0}

%

  \item  $\Sigma^{-1}\cX \subseteq \cX$ and $\Sigma\cY \subseteq \cY$.

  \item  $\cT( \cX,\cY ) = 0$.

  \item  For each object $t \in \cT$ there is a triangle
  $x \rightarrow t \rightarrow y$ with $x \in \cX$, $y \in \cY$. 

\end{enumerate}
The {\em co-heart} is $\cC = \cX \cap \Sigma^{-1}\cY$.

The co-t-structure is called {\em bounded} if
\[
\tag*{$\Box$}
  \bigcup_{i \in \BZ} \Sigma^i \cX = \bigcup_{i \in \BZ} \Sigma^i \cY = \cT.
\]
\end{Definition}

In contrast to t-structures, the objects $x$ and $y$ in Definition
\ref{def:co-t-structure}(iii) do {\em not} in general depend
functorially on $t$, see \cite[rmk.\ 1.2.2]{B}.

The following is the canonical example of a co-t-structure.

\begin{Example}
[The standard co-t-structure]
\label{exa:co-t-structure}
Let $R$ be a ring, $\cP = \Prj\, R$ the category of projective
left-modules over $R$.  Let $\cX$ and $\cY$ be the isomorphism
closures in the bounded homotopy category $\Kb( \cP )$ of the subsets 
\begin{equation}
\label{equ:co-t-structure_by_degree}
\vcenter{
  \xymatrix @-2pc @C=0.7pc @R=0.3pc
  {
    \{\, \cdots \ar[r] & 0 \ar[r] & 0 \ar[r] & 0 \ar[r] & P^0 \ar[r] & P^1 \ar[r] & P^2 \ar[r] & \cdots \,|\, P^i \in \cP, \, P^i = 0 \mbox{ for } i \gg 0 \,\} \lefteqn{,} \\
    \{\, \cdots \ar[r] & P^{-3} \ar[r] & P^{-2} \ar[r] & P^{-1} \ar[r] & 0 \ar[r] & 0 \ar[r] & 0 \ar[r] & \cdots \,|\, P^i \in \cP, \, P^i = 0 \mbox{ for } i \ll 0 \,\} \lefteqn{.} 
  }
         }
\end{equation}
We will show that if $\cX$ and $\cY$ are viewed as full subcategories,
then $( \cX,\cY )$ is a bounded co-t-structure with co-heart $\cC$
equivalent to $\cP$. 

Recall that $\cX$ and $\cY$ are required to be closed under
isomorphisms, direct sums, and direct summands.  The two former
properties are immediate, and the latter follows from the results in
\cite[secs.\ 3 and 4]{S}.

We next check the conditions in Definition \ref{def:co-t-structure}.
Conditions (i) and (ii) are clear.  The triangle in Definition
\ref{def:co-t-structure}(iii) can be obtained by hard truncation
of the object $t = \cdots \rightarrow P^{-2} \rightarrow P^{-1}
\rightarrow P^0 \rightarrow P^1 \rightarrow \cdots$ in $\Kb( \cP )$,
see Figure \ref{fig:2} in the introduction.

It is immediate that the co-t-structure $( \cX,\cY )$ is bounded.

Finally, it follows from \cite[cor.\ 4.11]{S} that the coheart $\cC =
\cX \cap \Sigma^{-1}\cY$ is equivalent to $\cP$.
\hfill $\Box$
\end{Example}

The term {\em silting set} was coined in \cite{KV}.  The following
definition was made in \cite[def.\ 2.1]{AI}.

\begin{Definition}
\label{def:silting}
A {\em silting subcategory} $\cC$ of $\cT$ is a full subcategory,
closed under isomorphisms, direct sums, and direct summands,
which satisfies
\begin{enumerate}
\setlength\itemsep{4pt}

  \item  $\cT( \cC,\Sigma^{ >0 }\cC ) = 0$.

  \item  Each object in $\cT$ can be obtained from $\cC$ by taking
    finitely many (de)suspensions, triangles, and direct summands.

\end{enumerate}
A {\em silting object} $s$ of $\cT$ is an object such that $\add( s )$
is a silting subcategory. 
\hfill $\Box$
\end{Definition}




The following was proved in \cite[cor.\ 5.9]{MSSS}.

\begin{Theorem}
[Mendoza Hern\'{a}ndez et.al.]
\label{thm:MSSS}
The map
\[
  ( \cX,\cY ) \mapsto \cC = \cX \cap \Sigma^{ -1 }\cY
\]
is a bijection between bounded co-t-structures and silting
subcategories of $\cT$.
\end{Theorem}

\begin{Remark}
The inverse map sends a silting subcategory $\cC$ to a pair $( \cX,\cY
)$ where $\cX$ is the smallest full subcategory, closed under
isomorphisms, direct sums, and direct summands, which is closed under
$\Sigma^{ -1 }$ and contains $\cC$.  Similarly, $\cY$ is the smallest
full subcategory, closed under isomorphisms, direct sums, and direct
summands, which is closed under $\Sigma$ and contains $\Sigma \cC$.
\hfill $\Box$
\end{Remark}

\begin{Example}
\label{exa:little_co-t-structure}
We continue Example \ref{exa:little_t-structure}, so
$\Lambda = \BC A_2$ is the path algebra of the quiver $A_2$ from
Equation \eqref{equ:A2} and $\cP = \prj\,\Lambda$ is the category of
finite dimensional projective left modules over $\Lambda$.

There is a bounded co-t-structure $( \cX,\cY )$ in $\Kb( \cP )$ where
$\cX$ and $\cY$ are the isomorphism closures in $\Kb( \cP )$ of the
subsets in Equation \eqref{equ:co-t-structure_by_degree}; this is
shown by the same method as in Example \ref{exa:co-t-structure}.

Recall that $\Lambda$ also has a bounded derived category $\Db( \cM
)$, see Example \ref{exa:little_t-structure}.  Since $\Lambda$ has
global dimension $1$, the triangulated categories $\Kb( \cP )$ and
$\Db( \cM )$ are equivalent, so $\Kb( \cP )$ has the Auslander--Reiten
quiver shown in Equation \eqref{equ:AR_quiver}.  We redraw the quiver,
this time with red and green vertices showing $\cX$ and $\cY$.
\begin{equation}
\label{equ:AR_quiver2}
\vcenter{
  \xymatrix @-1.0pc @! {
    & {\color{red} x_{ -3 }} \ar[dr] & & {\color{red} x_{ -1 }} \ar[dr] & & {\color{red} x_1} \ar[dr] & & {\color{green} x_3} \ar[dr] & & {\color{green} x_5} \ar[dr] & \\
    \cdots \ar[ur] & & {\color{red} x_{ -2 }} \ar[ur] & & {\color{red} x_0} \ar[ur] & & {\color{black} x_2} \ar[ur] & & {\color{green} x_4} \ar[ur] & & \cdots
               }
         }
\end{equation}
Note that $x_2$ is neither in $\cX$ nor $\cY$.  Indeed, if we abuse
notation to confuse isomorphism classes with individual objects, then
we can set $t = x_2$ and the triangle in Definition
\ref{def:co-t-structure}(iii) becomes $x_1 \rightarrow x_2 \rightarrow
x_3$.

The coheart $\cC = \cX \cap \Sigma^{ -1 }\cY$ is determined by
\[
  \cC = \add( x_0 \oplus x_1 ).
\]
Theorem \ref{thm:MSSS} implies that $\cC$ is a silting subcategory of
$\Kb( \cP )$.  The corresponding isomorphism class of silting objects
is $x_0 \oplus x_1$.  \hfill $\Box$
\end{Example}

\section{Towers which build an arbitrary object from objects of the (co-)heart}
\label{sec:towers}

The following two results were proved in \cite[p.\ 34]{BBD} and
\cite[prop.\ 1.5.6]{B}.  A wavy arrow $\xymatrix { s \ar@{~>}[r] & t
}$ denotes a morphism $s \rightarrow \Sigma t$.

\begin{Proposition}
[Beilinson, Bernstein, and Deligne]
\label{pro:t-tower}
Let $( \cA , \cB )$ be a bounded t-structure in $\cT$ with heart
$\cH$.  For each object $t \in \cT$, there is an integer $n \geqslant
1$ and a diagram consisting of triangles,
\[
  \xymatrix @-2.5pc @! {
    0 = t_0 \ar[rr] & & t_1 \ar[rr] \ar[dl] & & t_2 \ar[rr] \ar[dl] & & \cdots \ar[rr] & & t_{n-1} \ar[rr] & & t_n = t, \ar[dl] \\
    & \Sigma^{i_1}h_1 \ar@{~>}[ul] & & \Sigma^{i_2}h_2 \ar@{~>}[ul] & & & & & & \Sigma^{i_n}h_n \ar@{~>}[ul] & \\
               }
\]
with $h_m \in \cH$ for each $m$ and $i_1 > i_2 > \cdots > i_n$. 
\end{Proposition}

\begin{Proposition}
[Bondarko]
\label{pro:co-t-tower}
Let $( \cX , \cY )$ be a bounded co-t-structure in $\cT$ with co-heart
$\cC$.  For each object $t \in \cT$, there is an integer $n \geqslant
1$ and a diagram consisting of triangles,
\[
  \xymatrix @-2.5pc @! {
    0 = t_0 \ar[rr] & & t_1 \ar[rr] \ar[dl] & & t_2 \ar[rr] \ar[dl] & & \cdots \ar[rr] & & t_{n-1} \ar[rr] & & t_n = t, \ar[dl] \\
    & \Sigma^{i_1}c_1 \ar@{~>}[ul] & & \Sigma^{i_2}c_2 \ar@{~>}[ul] & & & & & & \Sigma^{i_n}c_n \ar@{~>}[ul] & \\
               }
\]
with $c_m \in \cC$ for each $m$ and $i_1 < i_2 < \cdots < i_n$. 
\end{Proposition}

\begin{Example}
\label{exa:t-tower}
Consider the t-structure in Example \ref{exa:t-structure}.  If $t \in
\Db( \cM )$ is given, then there is a diagram as in Proposition
\ref{pro:t-tower} where the objects $\Sigma^{ i_m }h_m$ are of the
form $\Sigma^i\H^{ -i }( t )$.  The diagram expresses that $t$ can be
built from its cohomology modules $\H^{ -i }( t )$.
\hfill $\Box$
\end{Example}

\begin{Example}
\label{exa:co-t-tower}
Consider the co-t-structure in Example \ref{exa:co-t-structure}.  If
$t = \cdots \rightarrow P^{-2} \rightarrow P^{-1} \rightarrow P^0
\rightarrow P^1 \rightarrow \cdots$ in $\Kb( \cP )$ is given, then
there is a diagram as in Proposition \ref{pro:co-t-tower} where the
objects $\Sigma^{ i_m }c_m$ are of the form $\Sigma^iP^{ -i }$.  The
diagram expresses that $t$ can be built from its constituent modules
$P^{ -i }$.  \hfill $\Box$
\end{Example}

\section{Categories skewed towards t- or co-t-structures}
\label{sec:skewed}

In this section, $d$ is a fixed integer.

\begin{Definition}
\label{def:spherical}
If $\cT$ is $\BC$-linear, then an object $s \in \cT$ is called {\em
  $d$-spherical} if there is an isomorphism
\[
  \cT( s,\Sigma^{ \ast }s ) \cong \BC[ X ] / ( X^2 )
\]
of graded algebras where $\deg X = d$.  
\hfill $\Box$
\end{Definition}

\begin{Remark}
If $d \neq 0$, then $s$ is $d$-spherical if and only if there are
isomorphisms of $\BC$-vector spaces
\[
\tag*{$\Box$}
  \cT( s,\Sigma^i s )
  \cong
  \left\{
    \begin{array}{cl}
      \BC & \mbox{ for } i \in \{ 0,d \}, \\[1mm]
      0   & \mbox{ otherwise. }
    \end{array}
  \right.
\]
\end{Remark}

The following was proved in \cite[thm.\ 2.1]{KYZ}.

\begin{Theorem}
[Keller, Yang, and Zhou]
There is a triangulated category $\cS_d$ which is algebraic,
$\BC$-linear with finite dimensional $\Hom$-spaces, and contains a
$d$-spherical object $s$ such that each object in $\cS_d$ can be
obtained from $s$ by taking finitely many (de)suspensions, triangles,
and direct summands.

Up to triangulated equivalence, $\cS_d$ is unique.
\end{Theorem}

The following was proved in \cite[thm.\ A]{HJY}.  It clearly implies
that if $d \leqslant 0$, then bounded co-t-structures are more useful
than bounded t-structures for the study of $\cS_d$.

\begin{Theorem}
[Holm, J, and Yang]
\label{thm:HJY}
If $d \leqslant 0$ then $\cS_d$ has no bounded t-structures.  It has one
familiy of bounded co-t-structures, all of which are (de)suspensions
of a canonical one.

If $d \geqslant 1$ then $\cS_d$ has no bounded co-t-structures.  It has
one familiy of bounded t-structures, all of which are (de)suspensions
of a canonical one.
\end{Theorem}

\section{The bijections of K\"{o}nig and Yang}
\label{sec:KY}

\begin{Definition}
A {\em simple minded collection} in $\cT$ is a set $\{ t_1, \ldots,
t_n \}$ of objects of $\cT$ with the following properties.
\begin{enumerate}
\setlength\itemsep{4pt}

  \item  $\cT( t_i,\Sigma^{ <0 }t_j ) = 0$ for all $i$ and $j$.

  \item  $\cT( t_i,t_i )$ is a division ring for each $i$ and $\cT(
    t_i,t_j ) = 0$ when $i \neq j$.  

  \item  Each object in $\cT$ can be obtained from $t_1$, $\ldots$,
    $t_n$ by taking finitely many (de)sus\-pen\-si\-ons, triangles,
    and direct summands.
    \hfill $\Box$

\end{enumerate}
\end{Definition}

Let $\Lambda$ be a finite dimensional $\BC$-algebra, $\cM = \mod\,
\Lambda$ the category of finite dimensional left modules over
$\Lambda$, and $\cP = \prj\, \Lambda$ the category of finite
dimensional projective left modules over $\Lambda$.  The following was
proved in \cite[thm.\ (6.1)]{KY}.  We interpret it as indicating that
the role of t-structures in derived categories is taken over by
co-t-structures in homotopy categories.

\begin{Theorem}
\label{thm:KY}
There are bijections between the following sets. 
\begin{enumerate}
\setlength\itemsep{4pt}

  \item  Bounded t-structures in $\Db( \cM )$ whose hearts are length
    categories (``length category'' means that each object has finite
    length). 

  \item  Bounded co-t-structures in $\Kb( \cP )$.

  \item  Isomorphism classes of simple minded collections in $\Db( \cM )$.

  \item  Isomorphism classes of basic silting objects in $\Kb( \cP )$
    (basic means no repeated indecomposable summands). 

\end{enumerate}
\end{Theorem}

\begin{Remark}
There is an extensive body of work on the bijections of Theorem
\ref{thm:KY} which predates \cite{KY}, see \cite{AI}, \cite{Al-N},
\cite{AST}, \cite{B}, \cite{BRT}, \cite{KN}, \cite{KV}, and
\cite{MSSS}.  The contributions of these papers to Theorem
\ref{thm:KY} are explained in the introduction to \cite{KY}.
\hfill $\Box$
\end{Remark}

\begin{Remark}
The proof of Theorem \ref{thm:KY} occupies a significant part of
\cite{KY}, and we only show how some of the bijections are defined.

(i) to (iii):  Let $( \cA,\cB )$ be a bounded t-structure in $\Db( \cM )$
whose heart $\cH = \cA \cap \Sigma \cB$ is a length category.  Take a
simple object from each isomorphism class of simple objects in $\cH$.
This gives a simple minded system in $\Db( \cM )$, see \cite[sec.\
5.3]{KY}.   

(ii) to (iv):  Let $( \cX,\cY )$ be a bounded co-t-structure in $\Kb(
\cP )$.  The co-heart $\cC = \cX \cap \Sigma^{ -1 }\cY$ is a silting
subcategory by Theorem \ref{thm:MSSS}, and there is a silting object
$s$ such that $\cC = \add( s )$, see \cite[sec.\ 5.2]{KY}.

(i) to (ii):  Let $( \cX,\cY )$ be a bounded co-t-structure in $\Kb(
\cP )$.  Set
\begin{align*}
  \cA = \{ a \in \Db( \cM )
        & \,|\, \Hom_{ \Db( \cM ) }( x,a ) = 0
          \mbox{ for each $x \in \cX$} \}, \\
  \cB = \{ b \in \Db( \cM )
        & \,|\, \Hom_{ \Db( \cM ) }( y,b ) = 0
          \mbox{ for each $y \in \cY$} \}
\end{align*}
and view these two sets as full subcategories of $\Db( \cM )$.  Then
$( \cA,\cB )$ is a bounded t-structure in $\Db( \cM )$ with length
heart, see \cite[sec.\ 5.7]{KY}.  
\hfill $\Box$
\end{Remark}

\section{The silting mutation of Aihara and Iyama}
\label{sec:mutation}

Silting mutation is an operation which changes one silting subcategory
into another.  By virtue of Theorem \ref{thm:MSSS}, it can be viewed
as changing one bounded co-t-structure into another.  This leads to
the definition of the so-called silting quiver of the triangulated
category $\cT$ which shows how silting subcategories, and hence
co-t-structures, fit together inside $\cT$.

In this section, $\cT$ is $\BC$-linear with finite dimensional
$\Hom$-spaces and split idempotents, and $m = m_0 \oplus m_1$ is a
basic silting object of $\cT$ with $m_0$ indecomposable.  The
following definition and theorem are special cases of \cite[sec.\
2.4]{AI}.

\begin{Definition}
[Aihara and Iyama]
\label{def:silting_mutation}
Let $r \stackrel{\rho}{\rightarrow} m_0
\stackrel{\lambda}{\rightarrow} \ell$ be a minimal right $\add( m_1
)$-approximation and a minimal left $\add( m_1 )$-approximation of
$m_0$.  Complete these morphisms to triangles
\[
  m_0^{\sim} \rightarrow r \stackrel{\rho}{\rightarrow} m_0
  \;\;,\;\;
  m_0 \stackrel{\lambda}{\rightarrow} \ell \rightarrow m_0^{\dagger}
\]
in $\cT$ and set
\[
  \mu^-( m,m_1 ) = m_0^{ \sim } \oplus m_1
  \;\;,\;\;
  \mu^+( m,m_1 ) = m_0^{ \dagger } \oplus m_1.
\]
These are called {\em right} and {\em left silting mutations of $m$}. 
\hfill $\Box$
\end{Definition}

\begin{Theorem}
[Aihara and Iyama]
\label{thm:silting_mutation}
\begin{enumerate}
\setlength\itemsep{4pt}

  \item  The silting mutations $\mu^-( m,m_1 )$ and $\mu^+( m,m_1 )$
    are basic silting objects of $\cT$.

  \item  Right and left silting mutations are inverse in the sense
    that
\[
  \mu^-\big( \mu^+( m,m_1 ),m_1 \big) \cong m
  \;\;,\;\;
  \mu^+\big( \mu^-( m,m_1 ),m_1 \big) \cong m.
\]

\end{enumerate}
\end{Theorem}

The following is a special case of \cite[def.\ 2.41]{AI}.

\begin{Definition}
[Aihara and Iyama]
\label{def:silting_quiver}
The {\em silting quiver} of $\cT$ has a vertex for each isomorphism
class of basic silting objects of $\cT$, and an arrow $[m] \rightarrow
[m^*]$ if $m^*$ is a left silting mutation of $m$, where square
brackets denote isomorphism class.
\hfill $\Box$
\end{Definition}

\begin{Remark}
The silting quiver gives a picture of how silting mutation moves from
one silting object to another, hence from one silting subcategory to
another.  By virtue of Theorem \ref{thm:MSSS}, it gives a picture of
how silting mutation moves from one bounded co-t-structure to another.
\hfill $\Box$
\end{Remark}

\begin{Example}
We continue Example \ref{exa:little_co-t-structure}, so $\Lambda = \BC
A_2$ is the path algebra of the quiver $A_2$ from Equation
\eqref{equ:A2}.  The bounded homotopy category $\Kb( \cP )$ has the
Auslander--Reiten quiver shown in Equation \eqref{equ:AR_quiver2}.

Recall from Example \ref{exa:little_co-t-structure} that $\Kb( \cP )$
has the isomorphism class of silting objects $x_0 \oplus x_1$.
Indeed, $x_0 \oplus x_1$ is a vertex in the silting quiver of $\Kb(
\cP )$.  The full quiver was determined in \cite[exa.\ 2.45]{AI}, see
Figure \ref{fig:AI}.
\begin{figure}
\[
  \xymatrix @!0 @C=2.5pc @R=3.0pc
  {
    \cdots \ar[dd] \ar[dddrrr] &&&&&& \cdots \ar[dddrrr] \ar[dddlll] |!{[ddll];[dd]}\hole &&&&&& \cdots \ar[dddlll] |!{[ddll];[dd]}\hole \\
    &&& x_{ -2 } \oplus x_2 \ar[dddrrr] \ar[dddlll] |!{[dl];[d]}\hole |!{[ddll];[dd]}\hole &&&&&& x_{ -5 } \oplus x_5 \ar[dddrrr] \ar[dddlll] |!{[dl];[d]}\hole |!{[ddll];[dd]}\hole &&&\\
    x_0 \oplus x_1 \ar[dd] \ar[dddrrr] &&&&&& x_{ -3 } \oplus x_4 \ar[dddrrr] \ar[dddlll] |!{[dl];[d]}\hole |!{[ddll];[dd]}\hole &&&&&& \cdots \ar[dddlll] |!{[dl];[d]}\hole |!{[ddll];[dd]}\hole \\
    &&& x_{ -1 } \oplus x_3 \ar[dddrrr] \ar[dddlll] |!{[dl];[d]}\hole |!{[ddll];[dd]}\hole &&&&&& x_{ -4 } \oplus x_6 \ar[dddrrr] \ar[dddlll] |!{[dl];[d]}\hole |!{[ddll];[dd]}\hole &&&\\
    x_1 \oplus x_2 \ar[dd] \ar[dddrrr] &&&&&& x_{ -2 } \oplus x_5 \ar[dddrrr] \ar[dddlll] |!{[dl];[d]}\hole |!{[ddll];[dd]}\hole &&&&&& \cdots \ar[dddlll] |!{[dl];[d]}\hole |!{[ddll];[dd]}\hole \\
    &&& x_{ 0 } \oplus x_4 \ar[dddrrr] \ar[dddlll] |!{[dl];[d]}\hole |!{[ddll];[dd]}\hole &&&&&& x_{ -3 } \oplus x_7 \ar[dddrrr] \ar[dddlll] |!{[dl];[d]}\hole |!{[ddll];[dd]}\hole &&&\\
    x_2 \oplus x_3 \ar[dd] \ar[dddrrr] &&&&&& x_{ -1 } \oplus x_6 \ar[dddrrr] \ar[dddlll] |!{[dl];[d]}\hole |!{[ddll];[dd]}\hole &&&&&& \cdots \ar[dddlll] |!{[dl];[d]}\hole |!{[ddll];[dd]}\hole \\
    &&& x_{ 1 } \oplus x_5 \ar[dddrrr] \ar[dddlll] |!{[dl];[d]}\hole &&&&&& x_{ -2 } \oplus x_8 \ar[dddrrr] \ar[dddlll] |!{[dl];[d]}\hole &&&\\
    x_3 \oplus x_4 \ar[dd] &&&&&& x_{ 0 } \oplus x_7 &&&&&& \cdots \\
    &&& x_{ 2 } \oplus x_6 &&&&&& x_{ -1 } \oplus x_9 &&&\\
    \cdots &&&&&& \cdots &&&&&& \cdots
  }
\]
  \caption{The silting quiver of the bounded homotopy category $\Kb(
    \cP )$ where $\cP$ is the category of finite dimensional projective
    modules over $\BC A_2$.} 
\label{fig:AI}
\end{figure}
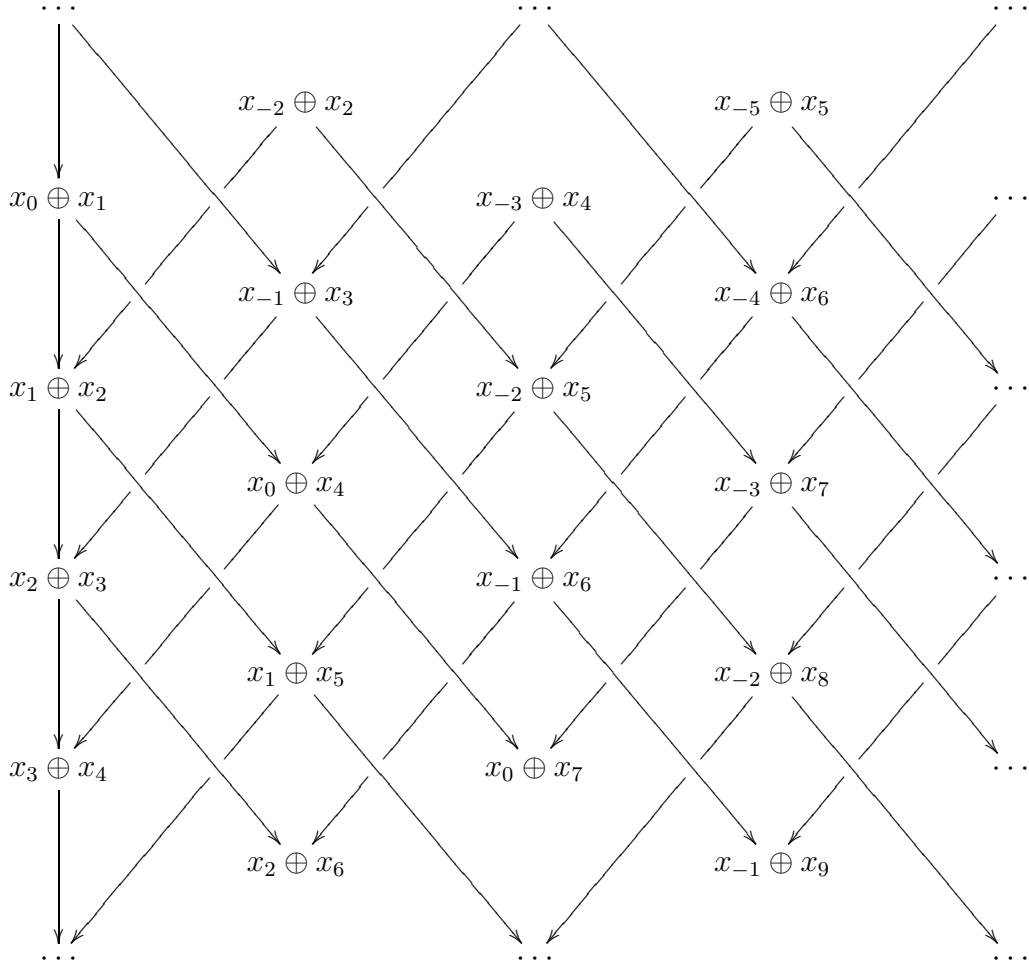
As the quiver shows, there is a left silting mutation of $x_0 \oplus
x_1$ which gives $x_0 \oplus x_4$.  The isomorphism classes of silting
objects $x_0 \oplus x_1$ and $x_0 \oplus x_4$ give rise to silting
subcategories $\add( x_0 \oplus x_1 )$ and $\add( x_0 \oplus x_4 )$
which, under the bijection of Theorem \ref{thm:MSSS}, correspond to
two bounded co-t-structures.  The first of these is shown on the AR
quiver of $\Kb( \cP )$ in Equation \eqref{equ:AR_quiver2}.  The second
co-t-structure $( \cX',\cY' )$ can be shown as follows, where the red
and green vertices show $\cX'$ and $\cY'$.
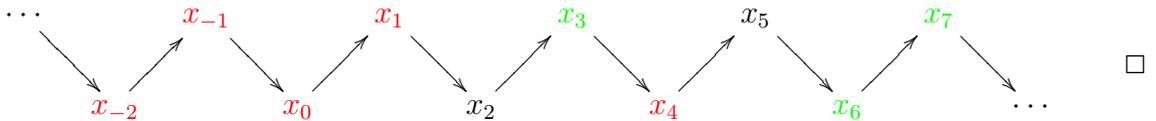
\[
\tag*{$\Box$}
\vcenter{
  \xymatrix @-1.2pc @! {
    \cdots \ar[dr] & & {\color{red} x_{ -1 }} \ar[dr] & & {\color{red}
      x_1} \ar[dr] & & {\color{green} x_3} \ar[dr] & &
      x_5 \ar[dr] & & {\color{green} x_7} \ar[dr] \\
    & {\color{red} x_{ -2 }} \ar[ur] & & {\color{red} x_0} \ar[ur] & &
    x_2 \ar[ur] & & {\color{red} x_4} \ar[ur] & &
    {\color{green} x_6} \ar[ur] & & \cdots
               }
         }
\]
\end{Example}

\medskip
\noindent
{\bf Acknowledgement.}
I thank Kiyoshi Igusa, Alex Martsinkovsky, and Gordana Todorov for
inviting me to speak at the Maurice Auslander International Conference
and to write this survey.

\end{document}